\documentclass{amsproc}

\usepackage{amssymb,latexsym,amsfonts,verbatim,amscd,ifthen}

\usepackage[usenames,dvipsnames]{pstricks}
\usepackage{pstricks}
\usepackage{epsfig}
\usepackage{pst-grad} 
\usepackage{pst-plot} 

\newtheorem{thm1}{Theorem}[section]
\newtheorem{lem1}[thm1]{Lemma}

\newtheorem{def1}[thm1]{Definition}

\newtheorem{cor1}[thm1]{Corollary}

\newtheorem{prop1}[thm1]{Proposition}

\newtheorem{example}[equation]{Example}

\DeclareMathOperator{\sgn}{sgn}
\DeclareMathOperator{\HH}{\tilde{H}}
\DeclareMathOperator{\supp}{supp}

\DeclareMathOperator{\Tor}{Tor}

\DeclareMathOperator{\id}{id}

 \DeclareMathOperator{\im}{Im}

 \def\mapright#1{\smash{\mathop{\longrightarrow}\limits^{#1}}}

\renewcommand\b{\beta}

\renewcommand\th{\theta}
\newcommand\D{\Delta}

\renewcommand{\L}{\mathcal{L}}
\newcommand{\A}{\mathcal{A}}

\newcommand{\om}{\omega}

\begin{document}

\title[Syzygies]
{  On simple $\A$-multigraded  minimal resolutions }
\author[H. Charalambous]{Hara Charalambous}
\address { Department of Mathematics, Aristotle University of Thessaloniki,
Thessaloniki \\54124, GREECE
 } \email{hara@math.auth.gr}
\author[A. Thoma]{ Apostolos Thoma}
\address { Department of Mathematics, University of Ioannina,
Ioannina 45110, GREECE } \email{athoma@cc.uoi.gr}

\keywords{Resolutions, lattice ideal, syzygies, indispensable
syzygies, Scarf complex.} \subjclass{13D02, 13D25}

\begin{abstract} Let $\A$ be a semigroup whose only invertible
element is 0.  For an $\A$-homogeneous ideal
 we discuss the notions of simple $i$-syzygies and simple minimal free resolutions
 of $R/I$.
When $I$ is a lattice ideal, the simple $0$-syzygies of $R/I$ are
the binomials in $I$.  We show that for an appropriate choice of
bases every $\A$-homogeneous minimal free resolution of $R/I$ is
simple. We introduce the gcd-complex $\D_{\gcd}(\bf b)$ for a
degree $\mathbf{b}\in \A$. We show that the homology of
$\D_{\gcd}(\bf b)$ determines the $i$-Betti numbers of degree $\bf
b$. We discuss the notion of an indispensable complex of $R/I$. We
show that the Koszul complex of a complete intersection lattice
ideal $I$ is the indispensable resolution of $R/I$ when the
$\A$-degrees of the elements of the generating $R$-sequence are
incomparable.

\end{abstract}
\maketitle

\date{}

\section{Notation}
\label{intro_section}

Let $ \mathcal{L}\subset \mathbb{Z}^n$ be a lattice  such that
$\L\cap \mathbb{N}^n=\{ {\bf 0}\}$ and let $\mathcal{A}$ be the
subsemigroup of  $\mathbb{Z}^n/\mathcal{L}$ generated by $\{{\bf
a_i}={\bf e}_i+\L:1\leq i\leq n\}$ where $\{e_i:\  1\leq i\leq
n\}$ is the canonical basis of $\mathbb{Z}^n $. Since the only
element in $\A$ with an inverse is 0, it follows that we can
partially order $\A$ with the relation
\[{\bf c} \geq {\bf d} \Longleftrightarrow \ \textrm{there is} \
{\bf e} \in \A \ \textrm{such that} \ {\bf c}={\bf d}+{\bf e}.\]

Let $\Bbbk$ be a field. We consider the polynomial ring
$R=\Bbbk[x_1,\ldots, x_n]$. We set $\deg_{\A}(x_i)={\bf a}_i$. If
${\bf x}^{{\bf v}}=x_1^{v_1} \cdots x_n^{v_n}$  then we set
\[ \deg_{ \A}({\bf x}^{{\bf v}}):=v_1{\bf
a}_1+\cdots+v_n{\bf a}_n \in \mathcal{A}\ .\] It follows that $R$
is positively multigraded by the semigroup $\mathcal{A}$, see
\cite{MS}. The {\em lattice ideal}  associated to $ \mathcal{L}$
is the ideal $I_{ \mathcal{L}}$ (or $I_\A$) generated by all the
binomials ${\bf x}^{{\bf u}_+}- {\bf x}^{{\bf u}_-}$ where
$\mathbf{u}_+, \mathbf{u}_-\in \mathbb{N}^n$ and
$\mathbf{u}=\mathbf{u}_+-\mathbf{u}_- \in \mathcal{L}$. We note
that if ${\bf x}^{{\bf u}_+}- {\bf x}^{{\bf u}_-}\in I_{
\mathcal{L}}$ then $\deg_{\A}{\bf x}^{{\bf u}_+}= \deg_{\A}{\bf
x}^{{\bf u}_-}$. Prime lattice ideals are the defining ideals of
toric varieties and are called toric ideals, \cite{St}. In general
lattice ideals arise in problems from diverse areas of
mathematics, including toric geometry, integer programming,
dynamical systems, graph theory, algebraic statistics,
hypergeometric differential equations, we refer to \cite{ES} for
more details.

We say than an ideal $I$ of $R$ is $\A$-{\it homogeneous} if it is
generated by $\A$-homogeneous polynomials, i.e.~polynomials whose
monomial terms have the same $\A$-degree. Lattice ideals are
clearly $\A$-{\it homogeneous}. For the rest of the paper $I$ is
an $\A$-{\it homogeneous} ideal. For $\bf{b} \in \A$ we let
$R[-{\bf b}]$ be the $\A$-graded free $R$-module of rank 1 whose
generator has $\A$-degree ${\bf b}$. Let
\[{(\bf F_\bullet,  \mathbf{\phi}}):
\quad 0\mapright{} F_p\mapright{\phi_p}\cdots
\cdots\mapright{}F_1\mapright{\phi_1} F_0\mapright{} R/I
\mapright{} 0,\] be a minimal  $\A$-graded free resolution of
$R/I$.   The $i$-Betti number of $R/I$ of $\A$-degree
$\mathbf{b}$, $\beta_{i,{\mathbf{b}}}(R/I)$, equals  the rank of
the $R$-summand of  $F_i$ of $\A$-degree $\mathbf{b}$:
\[\beta_{i,{\mathbf{b}}}(R/I)=\dim_\Bbbk\Tor_i(R/I,\Bbbk)_{\mathbf{b}}\]
and is an invariant of $I$, see \cite{MS}. The degrees ${\bf b}$
for which $\beta_{i,{\mathbf{b}}}(R/I)\neq 0$ are called {\em
i-Betti} degrees. The minimal elements of the set $\{{\bf b}:
\beta _{i, {\bf b}}(R/I)\not= 0\}$ are called minimal {\em i-Betti
degrees}. The elements of $\im \phi_{i+1}=\ker \phi_i$ are the
$i$-syzygies of $R/I$ in $\bf F_\bullet$.

The problem of obtaining an explicit minimal free resolution of
$R/I$ is extremely difficult. One of the factors that make this
problem hard to attack, is that given a minimal free resolution
one can obtain by a change of basis a different description of
this resolution. To obtain some control over this, in \cite{ChTh}
we defined   {\it simple } minimal free resolutions. We also
defined and studied the gcd-complex $\D_{\gcd}(\bf b)$ for a
degree $\mathbf{b}\in \A$. We used this complex to generalize the
results in \cite{PS1} and to  construct the generalized algebraic
Scarf complex based on the connected components of $\D_{\gcd}(\bf
b)$ for degrees $\mathbf{b}\in \A$. When $I$ is a lattice ideal we
showed that the generalized algebraic Scarf complex is present in
every simple minimal free resolution of $R/I$. This current paper
analyzes in more detail the notions presented in \cite{ChTh}. We
note that the original motivation for this work came from a
question in Algebraic Statistics concerning conditions for  the
uniqueness of a minimal binomial generating set of toric ideals.

 The structure of this paper is as follows. In section
\ref{simple_syzygies} we discuss the notion of {\it simple}
$i$-syzygies of $R/I$. The simple $0$-syzygies of $R/I$ when $I$
is a lattice ideal are exactly the binomials of $I$. We also
discuss the notion of a simple minimal free resolution of $R/I$.
This notion requires the presence of a system of bases for the
free modules of the resolution. We show that for an appropriate
choice of bases every $\A$-homogeneous minimal free resolution of
$R/I$ is simple. In section \ref{gcd_complex} we discuss the
gcd-complex $\D_{\gcd}(\bf b)$ for a degree $\mathbf{b}\in \A$. We
show that the homology of $\D_{\gcd}(\bf b)$ determines the
$i$-Betti numbers of degree $\bf b$.  We count the numbers of
binomials that could be part of a minimal binomial generating set
of a lattice ideal up to a constant multiple. In section
\ref{indispensable_syzygies} we discuss the notion of {\it
indispensable $i$-syzygies}. Intrinsically indispensable
$i$-syzygies are present in all $\A$-homogeneous simple minimal
free resolutions. For the $0$-step and for a lattice ideal $I_\L$
this means that there are some binomials of the ideal $I_\L$ that
are part (up to a constant multiple) of all $\A$-homogeneous
systems of minimal binomial generators of $I_\L$.  A {\it strongly
indispensable $i$-syzygy} needs to be present in every minimal
free resolution of $R/I$ even if the resolution is not simple. For
the $0$-step and for a lattice ideal $I_\L$ this means that there
are some elements of  $I_\L$ that are part (up to a constant
multiple) of all $\A$-homogeneous minimal sets of generators of
$I_\L$, where the generators are not necessarily binomials. We
consider conditions for strongly indispensable $i$-syzygies to
exist. We show that the Koszul complex of a complete intersection
lattice ideal $I$ is indispensable when the $\A$-degrees of the
elements of the generating $R$-sequence are incomparable.

\section{Simple syzygies}\label{simple_syzygies}

We recall and generalize the definition of  a simple $i$-syzygy,
see \cite[Definition 3.1]{ChTh}, to arbitrary elements of an
$\A$-graded free module. Let $F$ be a free $\A$-graded module of
rank $\b$ and let $B=\{ E_{t}:\ t=1,\ldots, \b\}$ be an
$\A$-homogeneous basis of $F$. Let $h$ be an $\A$-homogeneous
element of $F$:
\[
h=\sum_{1\le t\le \b} (\sum_{ c_{{\bf a_t}\neq 0}} c_{\bf a_t}
{\bf x}^{\bf a_t})E_{i}\ .
\]
The $S$-support of $h$ with respect to $B$ is the set \[
S_{B}(h)=\{ {\bf x}^{\bf a_t}E_{i}:\ c_{\bf a_t}\neq 0 \}\ .\] We
introduce a partial order on the elements of $F$:
\[h'\le h  \textrm{ if and only if } S_{B}(h')\subset S_{B}(h)\ .\]

\begin{def1}{\rm  Let $F$ and $B$ be as above,
let $G$ be an $\A$-graded subset of $F$ and let $h$ be an
$\A$-homogeneous nonzero element of $ G$. We say that $h$ is
{\emph simple} in $G$ with respect to $B$ if  there is no nonzero
$\A$-homogeneous $ h'\in G$ such that $h'<h$.}
\end{def1}

In \cite[Theorem 3.4]{ChTh} we showed that if $(\bf F_\bullet, \phi)$ is a
minimal free resolution of $R/I$ then for any given basis $B$ of
$F_i$ there exists a minimal $\A$-homogeneous generating set of
$\ker\phi_i$ consisting of simple $i$-syzygies with respect to
$B$. The proof of the next proposition is an immediate
generalization of the proof of  that theorem and is omitted.

\begin{prop1}\label{existence_simple}
Let $F$ be a  free $\A$-graded module,   let $B$ be an
$\A$-homogeneous basis of $F$ and let $G$ be  an $\A$-graded
submodule of $F$. There is a minimal system of generators of $G$
each being simple in $G$ with respect to $B$.
\end{prop1}

Given an $\A$-homogeneous complex  of free modules $(\bf
G_\bullet, \phi)$ we specify $\A$-homoge\-neous bases $B_i$ for
the homological summands $G_i$. The collection of theses bases
forms a {\it system of bases} $\bf B$. We write ${\bf B}=(B_i)$
and we say that $B_i$ is in $\bf B$.

\begin{def1}{\rm A {\it based complex} $(\bf G_\bullet,
\phi, B)$ is an $\A$-homogeneous complex $(\bf G_\bullet, \phi)$
together with a system  of bases ${\bf B}=(B_i)$. Let $(\bf
G_\bullet, \phi, B)$ and $(\bf F_\bullet, \phi, C)$ be two based
complexes, ${\bf B}=(B_i)$ and ${\bf C}=(C_i)$. We say that the
complex homomorphism $\om: \bf G_\bullet\mapright{} \bf F_\bullet$
is a {\it based homomorphism} if  for each $E\in B_i$, there
exists an $H\in C_i$ such that $\om(E)=cH$ for some $c\in
\Bbbk^*$.}
\end{def1}

Let $I$ be an $\A$-homogeneous ideal and let $({\bf
F_\bullet},\phi)$ be a minimal $\A$-graded free resolution of
$R/I$. We let $s$ be the projective dimension of $R/I$ and $\b_i$
be the rank of $F_i$. For each $i$ we suppose that $B_i$  is an
$\A$-homogeneous basis of $F_i$ and we let $\mathbf{B}=(B_0,
B_1,\ldots, B_s)$.

\begin{def1}{\rm (\cite[Definition 3.5]{ChTh}) Let $({\bf
F_\bullet,\phi,B})$ be as above. We say that
 $({\bf
F_{\bullet},\phi, B})$ is {\it simple}  if and only if for each
$i$ and each $E\in B_i$, $\phi_i(E)$ is simple in $\ker\phi_{i-1}$
with respect to $B_{i-1}$.}
\end{def1}

We remark that when $I$ is a lattice ideal then for any choice of
basis $B_0$, the simple $0$-syzygies of $R/I$ are  the binomials
of $I$. It is an immediate consequence of Proposition
\ref{existence_simple} that one can construct a minimal simple
resolution of $R/I$ with respect to a system ${\bf B}
=(B_0,\ldots, B_s)$ starting with $B_0=\{1\}$, see also
\cite[Corollary 3.6]{ChTh}. In the next proposition we show that
any minimal free resolution $(\bf F_\bullet, \phi)$ of $R/I$
becomes simple with the right choice of bases.

\begin{prop1}\label{exist_simple_basis} Let $I$ be an $\A$-homogeneous ideal and
let $(\bf F_\bullet, \phi)$ be a minimal free resolution of $R/I$.
There exists a system of bases $\bf B$ so that $(\bf F_\bullet,
\phi,B)$ is simple.
\end{prop1}

\begin{proof} Let $C_0=\{1\}$ and for each $i>0$ choose a basis
$C_i=\{ H_{ti}: \ t=1,\ldots, \b_i\}$ of $F_i$. Let $(\bf
G_\bullet, \th)$ be a simple minimal free resolution of $R/I$ with
respect to ${\bf D}=(D_0, \ldots, D_s)$ where $D_0=\{1\}$ and
$D_i=\{E_{ti} :\ t=1,\ldots, \b_i\}$. Since $\bf G_\bullet$, $\bf
F_\bullet$ are both minimal projective resolutions of $R/I$ there
is an isomorphism of complexes $\bf h_\bullet:\ \bf
G_\bullet\mapright{} \bf F_\bullet$ that extends the identity map
on $R/I$. In particular $h_0=id_R$. For each $i$ we let
$H'_{ti}=h_i(E_{ti})$ and consider the set $B_i=\{ H'_{ti}: \
t=1,\ldots, \b_i\}$. We note that $B_0=\{1\}$. It is immediate
that $B_i$ is a basis for $F_i$. We claim that $({\bf F_\bullet,
\phi, B})$ is simple.

Indeed for $t=1,\ldots, \b_1$ using the commutativity of the
diagram we get that
\[\phi_1(H'_{t1})=\phi_1(h_1(E_{t1}))=h_0(\th_1(E_{t1}))=\th_1(E_{t1})\ .\]
Since $\th_1(E_{t1})$ is simple with respect to $C_0$ it follows
at once that  $\phi_1(H'_{t1})$ is simple with respect to $B_0$.
For $i>1$ and $t=1,\ldots, \b_i$ we have that
\[\phi_i(H'_{ti})=\phi_i(h_{i}(E_{ti}))=h_{i-1}(\th_i(E_{ti}))\ .\]
Suppose that $\phi_i(H'_{ti})$ were not simple with respect to
$B_{i-1}$. Since $h_{i-1}$ is  bijective   it follows that
$\th_i(E_{ti})$ is not  simple with respect to $D_{i-1}$, a
contradiction.
\end{proof}

Let $I$ be an   $\A$-homogeneous ideal and let $(\bf F_\bullet,
\phi)$ be a minimal free resolution of $R/I$. An $i$-syzygy $h$ of
$R/I$  {\it minimal} if $h$ is part of a minimal generating set of
$\ker\phi_i$. By the graded version of Nakayama's lemma it follows
that $h$ is {minimal} if and only if $h$ cannot be written as an
$R$-linear combination of $i$-syzygies of $R/I$ of strictly
smaller $\A$-degrees.

The next theorem examines the cardinality of the set of minimal
$i$-syzygies of a free resolution of $R/I$.

\begin{thm1}\label{finiteness_equiv_classes} Let $I$ be an
$\A$-homogeneous ideal and let
 $({\bf F_\bullet,
\phi})$ be a minimal free resolution of $R/I$. Let $\equiv$ be the
following equivalence relation among the elements of $F_i$:
 \[ h\equiv h' \text{ if and only if } h=ch',
c\in \Bbbk^*\ ,\] and let $B_i$ be a basis of $F_i$. The set of
equivalence classes of the $i$-syzygies of $R/I$ that are minimal
and simple with respect to $B_{i}$ is finite.
\end{thm1}

\begin{proof} We will show that the number of equivalence classes of the
$i$-syzygies that are simple and have $\A$-degree equal to an
$(i+1)$-Betti degree $\bf b$ of $R/I$ is finite.
 By \cite[Theorem
3.8]{ChTh} if $h, h'\in \ker \phi_i$ are simple with respect to
$B_i$ and $S_{B_i}(h)=S_{B_i}(h')$ then $h\equiv h'$. Thus it is
enough to show that there is only a finite number of candidates
for $S_{B_i}(h)$ when $h\in F_i$ has $\deg_A(h)=\bf b$.
 We consider the set
\[ C=\{ x^{\bf a}E_{t}: \deg_\A(x^{\bf a}E_{t})={\bf b} \ , \
E_{t}\in B_i\}.\] We note that $S_{B_i}(h)\subset \mathcal{P}(C)$
where $\mathcal{P}(C)$ is the power set of $C$. The number of
basis elements $E_{t}\in B_i$ such that $\deg_\A (E_{t})\le {\bf
b}$ is finite. Moreover for such $E_{t}$ the number of monomials $
x^{\bf a}$ such that $\deg_\A(x^{\bf a})+\deg_A( E_{t})={\bf b}$
is finite. It follows that $C$ and its power set $\mathcal{P}(C)$
are finite as desired.
 \end{proof}

In the next section we determine the cardinality of this set when
$I$ is  a lattice ideal and $i=0$. In other words we compute the
number of the equivalence classes of the minimal binomials of $I$.

\section{The gcd-complex}\label{gcd_complex}

For $\bf{b} \in \mathcal{A}$, we let $C_{\bf b}$ equal the fiber
\[C_{\bf
b}:=\deg_{\A}^{-1}(\bf b)=\deg_\mathcal{L}^{-1}({\bf b}):=\{x^u:\;
\deg_{\A}(x^u)=\bf b\}.\] Let $I_\L$ be a lattice ideal. The fiber
$C_{\bf b}$ plays an essential role in the study of the minimal
free resolution of $R/I_\L$ as is evident from several works, see
\cite{BaSt, ChKT, ChTh, DS, PS1, PS2}.

We denote the support of the vector ${\bf u}=(u_j)$ by $\supp({\bf
u}):= \{ i: u_i\neq 0\}$. Next we recall the definition of the
simplicial complex $\D_{\bf b}$ on $n$ vertices, constructed from
$C_{\bf b}$ as follows:
\[\D_{\bf b}:=\{F\subset\supp(\bf a):\ x^{\bf a}\in C_{\bf b}\}.\] $\D_{\bf b}$ has been studied
extensively, see for example \cite{AH, BCMP2, BH, CG, CP, P, PVT}.
Its homology determines the Betti numbers of $R/I_\L$:
\[ \b_{i,{\mathbf{b}}}(R/I_\L)=\dim_\Bbbk\HH_{i}(\D_{\bf b})\ ,\]
see \cite{Stanley} or \cite{MS} for a proof.

In this section we present another simplicial complex, the
gcd-complex $\D_{\gcd}(\bf b)$, whose construction is based upon
the divisibility properties of the monomials of $C_{\bf b}$.

\begin{def1}{\rm For a vector $\bf b\in \A$ we define  the {\it gcd-complex}
$\D_{\gcd}(\bf b)$ to be the simplicial complex with vertices the
elements of the fiber $C_{\bf b}$ and faces all subsets $T\subset
C_{\bf b} $ such that $\gcd(x^{\bf a}:\ x^{\bf a}\in T)\neq 1$.}
\end{def1}

The example below compares graphically the two simplicial
complexes in a particular case.

\begin{example}{\rm Let $R=\Bbbk[a,b,c,d]$ and let $\A$ be the semigroup
 $\A=\langle (4,0),(3,1),\\ (1,3), (0,4)\rangle$.
For ${\bf b}=(6,10)$ we consider the fiber $C_{(6,10)}= \{ bc^3,
ac^2d, b^2d^2\}$ and the corresponding simplicial complexes. We
see that

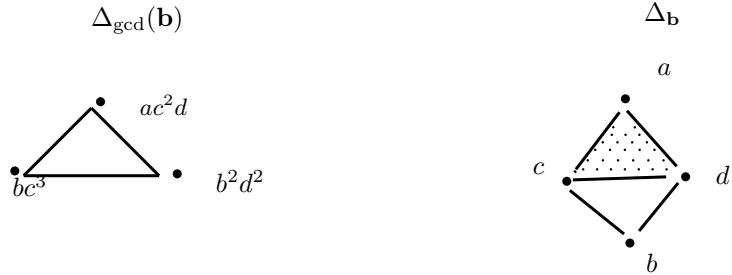
\scalebox{1} 
{
\begin{pspicture}(0,-1.8429687)(10.002812,1.8429687)
\usefont{T1}{ptm}{m}{n}
\rput(1.9123437,1.5745312){$\D_{\gcd}(\bf b)$}
\usefont{T1}{ptm}{m}{n}
\rput(8.8923435,1.6545312){$\D_{\bf b}$}
\psdots[dotsize=0.12](0.2809375,-0.43546876)
\psdots[dotsize=0.12](2.4409375,-0.47546875)
\psdots[dotsize=0.12](1.4209375,0.48453125)
\psline[linewidth=0.04cm](0.4,-0.5)(2.2,-0.5)
\psline[linewidth=0.04cm](0.4,-0.5)(1.3,0.4)
\psline[linewidth=0.04cm](1.3,0.4)(2.2,-0.5)
\usefont{T1}{ptm}{m}{n}
\rput(2.2584374,0.44453126){\small $ac^2d$}
\usefont{T1}{ptm}{m}{n}
\rput(0.48234376,-0.62546873){$bc^3$}
\usefont{T1}{ptm}{m}{n}
\rput(3.2823439,-0.5654687){$b^2d^2$}
\psdots[dotsize=0.12](8.400937,0.52453125)
\psdots[dotsize=0.12](7.6209373,-0.5754688)
\psdots[dotsize=0.12](9.200937,-0.5154688)
\psdots[dotsize=0.12](8.4609375,-1.3954687)
\psline[linewidth=0.04cm](8.320937,0.34453124)(7.7409377,-0.41546875)
\psline[linewidth=0.04cm](7.7409377,-0.43546876)(7.7209377,-0.45546874)
\psline[linewidth=0.04cm](8.420938,0.36453125)(9.080937,-0.37546876)
\psline[linewidth=0.04cm](7.6609373,-0.6954687)(8.380938,-1.2754687)
\psline[linewidth=0.04cm](9.100938,-0.59546876)(8.580937,-1.2354687)
\psline[linewidth=0.04cm](7.7009373,-0.55546874)(8.940937,-0.5154688)
\psline[linewidth=0.034cm,linestyle=dotted,dotsep=0.16cm](8.220938,0.14453125)(8.540937,0.14453125)
\psline[linewidth=0.034cm,linestyle=dotted,dotsep=0.16cm](8.200937,0.04453125)(8.600938,0.04453125)
\psline[linewidth=0.034cm,linestyle=dotted,dotsep=0.16cm](8.080937,-0.05546875)(8.700937,-0.05546875)
\psline[linewidth=0.034cm,linestyle=dotted,dotsep=0.16cm](7.9809375,-0.15546875)(8.800938,-0.15546875)
\psline[linewidth=0.034cm,linestyle=dotted,dotsep=0.16cm](7.8809376,-0.29546875)(8.920938,-0.29546875)
\psline[linewidth=0.034cm,linestyle=dotted,dotsep=0.16cm](7.8209376,-0.41546875)(8.9609375,-0.41546875)
\usefont{T1}{ptm}{m}{n}
\rput(8.922344,0.91453123){$a$}
\usefont{T1}{ptm}{m}{n}
\rput(7.242344,-0.40546876){$c$}
\usefont{T1}{ptm}{m}{n}
\rput(9.692344,-0.48546875){$d$}
\usefont{T1}{ptm}{m}{n}
\rput(8.752344,-1.6654687){$b$}
\end{pspicture}
} }
\end{example}

The main theorem of this section, Theorem \ref{iso_simp_compl} was
proved independently in \cite{OjVi}.

\begin{thm1}\label{iso_simp_compl} Let $\bf b\in \A$. The gcd complex
$\D_{\gcd}({\bf b})$ and the complex $\D_{\bf b}$ have the same
homology.
\end{thm1}
\begin{proof}
First we consider the simplicial complex $\D$ with vertices the
elements of the set $S=\{ \supp(\bf{a}): \ x^{\bf a}\in C_{\bf
b}\}$ and faces all subsets $T\subset S$ such that
\[\bigcap_{\supp(\bf{a})\in T} \supp(\bf{a})\neq \emptyset\ .\]
We define an equivalence relation among the vertices of
$\D_{\gcd}(\bf b)$: we let $x^{\bf{a}} \equiv x^{\bf{a'}}$ if and
only if $\supp(\bf{a})=\supp(\bf{a'})$. We note that the
subcomplex $A$ of $\D_{\gcd}(\bf b)$ on the vertices of an
equivalence class is contractible. By the Contractible Subcomplex
Lemma \cite{BW} we get that the quotient map $\pi: |\D_{\gcd}(\bf
b)|\mapright{} |\D_{\gcd}(\bf b)|/|A|$ is a homotopy equivalence.
A repeated application of the Contractible Subcomplex Lemma yields
that $\D_{\gcd}(\bf b)$ and $\D$ have the same homology.

Next we consider the family ${\mathcal F}$ of the  facets of
$\D_{\bf b}$ and the corresponding nerve complex $N({\mathcal
F})$. The vertices of $N({\mathcal F})$  correspond to the facets
of $\D_{\bf b}$, while the faces of $N({\mathcal F})$ correspond
to  collections of facets with nonempty intersection. It follows
that $N({\mathcal F})$ is isomorphic to $\D$. By \cite[Theorem
7.26]{Ro} the two complexes $\D_{\bf b}$ and $\D$ have the same
homology and the theorem now follows.
\end{proof}

The following is now immediate:

\begin{cor1}\label{betti_gcd} Let $I_\L$ be a lattice ideal.
\[ \b_{i,{\mathbf{b}}}(R/I_\L)=\dim_\Bbbk\HH_{i}(\D_{\gcd}({\bf b})).\]
\end{cor1}

The connected components of $\D_{\gcd}({\bf b})$ were used in
\cite{ChTh} to determine certain  complexes associated to a simple
minimal free $\A$-homogeous resolution of $R/I_\L$, see
\cite[Definitions 4.7 and 5.1]{ChTh}.  In Theorem
\ref{card_binomials} below we  use the complex $\D_{\gcd}({\bf
b})$ to determine the number of equivalence classes of minimal
binomial generators of $I_\L$. First we prove the following lemma:

\begin{lem1}\label{comp_graphs} For ${\bf b}\in \A$,   let $I_{\L, {\bf b}}$ be the ideal
generated by all binomials of $I_\L$ of $\A$-degree strictly
smaller than ${\bf b}$.  Let $G({\bf b})$ be the graph with
vertices the elements of $C_{\bf b}$   and edges all the sets
$\{x^{\bf u}, x^{\bf v}\}$ whenever $x^{\bf u}-x^{\bf v}\in I_{\L,
{\bf b}}$. A set of monomials in $C_{\bf b}$ forms the vertex set
of a component of  $G({\bf b})$ if and only if it forms the vertex
set of a component of $\D_{\gcd}({\bf b})$.
\end{lem1}

\begin{proof} We note that if
$x^{\bf u}, x^{\bf v}$ belong to the  same component of
$\D_{\gcd}({\bf b})$ then there exists a sequence of monomials
$x^{\bf u}=x^{{\bf u}_1}, x^{{\bf u}_2},\dots, x^{{\bf
u}_s}=x^{\bf v}$ such that $d=\gcd(x^{{\bf u}_i}, x^{{\bf
u}_{i+1}})\neq 1$. Therefore \[x^{{\bf u}_i}- x^{{\bf u}_{i+1}}=d(
\frac{x^{{\bf u}_i}}{d}- \frac{x^{{\bf u}_{i+1}}}{d}) \in I_{\L,
{\bf b}}\ .\] It follows that $x^{\bf u}-x^{\bf v}\in I_{\L, {\bf
b}}$ and $x^{\bf u}, x^{\bf v}$ belong to the  same component of
$G({\bf b})$.

For the converse we note that the binomials of degree $\bf b$ in
$I_{\L, {\bf b}}$ are spanned by binomials of the form $x^{\bf
a}(x^{\bf r}-x^{\bf s})$ where $x^{\bf a}\neq 1$. Moreover any
such binomial determines an edge from $x^{\bf a+r}$ to $x^{\bf
a+s}$ in $\D_{\gcd}({\bf b})$. Thus if $x^{\bf u}$, $x^{\bf v}$
lie in the same component of $G({\bf b})$ then  any minimal
expression of $x^{\bf u}-x^{\bf v}$ as a sum of binomials $ x^{\bf
a}(x^{\bf r}-x^{\bf s})$  results in a path from $x^{\bf u}$ to $
x^{\bf v}$ in $\D_{\gcd}({\bf b})$.
\end{proof}

The graph $G({\bf b})$ was first introduced in \cite{ChKT} to
determine the number of different binomial generating sets of a
toric ideal $I_\L$. The results stated for toric ideals in
\cite{ChKT}  hold more generally for lattice ideals with identical
proofs. We choose an ordering of the connected components of
$\D_{\gcd}({\bf b})$ and let $t_i({\bf b})$ be the number of
vertices of the $i$-th component of $\D_{\gcd}({\bf b})$.

\begin{thm1}\label{card_binomials} Let $I_\L$ be a lattice ideal and consider the
equivalence relation on $R$ of Theorem
\ref{finiteness_equiv_classes}. The cardinality of the set $T$ of
equivalence classes of the minimal binomials of $I_\L$ is given by
\[ |T|= \sum_{{\bf b}\in \A}\sum_{i\neq j}t_i({\bf b})t_j({\bf b}). \]
\end{thm1}

\begin{proof} In the course of the proof of
\cite[Theorem 2.6]{ChKT} applied to the lattice ideal  $I_\L$  it
was shown that the minimal binomials of $\A$-degree ${\bf b}$ are
the difference of monomials that belong to different connected
components of $G({\bf b})$. Lemma \ref{comp_graphs} and a counting
argument finishes the proof.
\end{proof}

We remark that if ${\bf b}\in \A$ is not a $1$-Betti degree of
$R/I_\L$, then there is no minimal binomial generator of
$\A$-degree ${\bf b}$. It follows that  $\D_{\gcd}({\bf b})$ has
exactly one connected component. The nontrivial contributions to
the formula of Theorem \ref{card_binomials} come from the
$1$-Betti degrees of $R/I_\L$.

\section{Indispensable syzygies}\label{indispensable_syzygies}

In this section we discuss the notion of indispensable complexes
that first appeared in \cite[Definition 3.9]{ChTh}. Intrinsically
an indispensable complex of $R/I$ is a based complex $
\bf(F_\bullet, \phi,B)$ that is ``contained" in any based simple
minimal free resolution of $R/I$.

The indispensable binomials  of a lattice ideal $I_\L$ are the
binomials that appear in every minimal system of binomial
generators of the ideal up to a constant multiple. They were first
defined in \cite{H-O} and their study was originally motivated
from Algebraic Statistics; see \cite{ATY, Hi-O, H-O, AT} for a
series of related papers.

\begin{thm1}\label{ind_0_syz} Let $I_\L$ be a lattice ideal. The indispensable binomials
of $I_\L$ occur exactly in the minimal $\A$-degrees $\bf b$ such
that $\D_{\gcd}(\bf b)$ consists of two disconnected vertices.
\end{thm1}

\begin{proof} This theorem was proved in  \cite{ChKT} for toric
ideals. The same proof applies to lattice ideals.
\end{proof}

An immediate consequence of Theorem \ref{ind_0_syz}  is the
following:

\begin{cor1}\label{start_indispens_complex} Let $I_\L$ be a lattice ideal and $S$
a minimal system of $\A$-homogeneous (not necessarily binomial)
generators  of $I_\L$. If $f$ is an indispensable binomial of
$I_\L$ then there is a $c\in \Bbbk^*$ such that $cf\in S$.
\end{cor1}

\begin{proof} Let $f$ be an indispensable binomial of
$I_\L$ and ${\bf b}=\deg_A f$. Since  $H_1(\D_{\gcd}({\bf b}))=1$
there is a unique element $f'$ in $S$ of $\A$-degree $\bf b$.
Since $C_{\bf b}$ is a set with exactly two elements it follows
that $f' $ is a binomial. Since $I_\L$ contains no monomials,  it
follows that $f'=cf$ for some $c\in \Bbbk^*$.
\end{proof}

It is clear that if $\bf(F_\bullet,\phi, B)$ is a minimal free resolution
of $R/I_\L$ and $f$ is an indispensable binomial, then there
exists an element $E\in B_1$ and a $c\in \Bbbk^*$ such that
$\phi_1(E)=cf$. We let the {\it indispensable $0$-syzygies} of
$R/I_\L$ to be the indispensable binomials of $I_\L$. We extend
the definition for $i\ge 0$:

\begin{def1}\label{indispensable_complex_def}{\rm  Let $\bf(F_\bullet, \phi,B)$ be a based complex.
We say that $(\bf F_\bullet, \phi, B)$  is an {\it indispensable
complex} of $R/I$ if
 for each based minimal simple  free resolution $\bf (G_\bullet,\theta,C)$ of
       $R/I$ where $C_0=\{ 1\}$, there is an injective based homomorphism $\bf \omega: \bf(F_\bullet, \phi,B) \to
       \bf (G_\bullet,\theta,C)$
       such that $\omega_0 = id_R$.
If ${\bf B}=(B_j)$ and $E\in B_{i+1}$ we say that
$\phi_{i+1}(E)\in F_i$ is an  \emph{indispensable i-syzygy} of
$R/I$.}
 \end{def1}

It follows immediately from the definition that an indispensable
$i$-syzygy of $R/I$ is simple.  Moreover if $(\bf F_\bullet, \phi,
B)$ is an {\it indispensable complex} of $R/I$ and $\bf
(G_\bullet,\theta,W)$ is a minimal simple free resolution of $R/I$ then
the based homomorphism of Definition
\ref{indispensable_complex_def} is unique, up to rearrangement of
the bases elements of the same $A$-degree and constant factors. In \cite[Theorem
5.2]{ChTh} we showed that if $I_\L$ is a lattice ideal then the
generalized algebraic Scarf complex is an indispensable complex.

The next theorem examines when the Koszul complex of a lattice
ideal generated by an $R$-sequence of binomials is indispensable.
Let  $I$ be an ideal generated by an $R$-sequence $f_1,\ldots,
f_s$ and let $(\bf K_\bullet,\bf \phi)$ be the Koszul complex on
the $f_i$. We denote the basis element $e_{j_1}\wedge \cdots\wedge
e_{j_t}$ of $K_t$ by $e_J$ where $J$ is the ordered set
$\{j_1,\ldots, j_t\}$ and let $\sgn[j_k,J]=(-1)^{k+1}$. For each
$j\in J$ we write $J_j$ for the set $J\setminus \{j\}$. The
canonical system of bases ${\bf B}=(B_0, \ldots, B_s)$ consists of
the following: $B_0=\{ 1\}$, $B_1=\{ e_i: i=1,\ldots, s\}$ where
$\phi_1(e_i)=f_i$ and $B_t=\{e_{J}:\ J=\{j_1,\ldots,j_t\},\  1 \le
j_1< \ldots < j_t\le s\}$ where
\[ \phi_{t}(e_{J})= \sum_{j\in J} \sgn[j,J] \ f_{j}
e_{J_j}\ .\] In \cite[Example 3.7]{ChTh} it was shown that $(\bf
K_\bullet,\bf \phi,\bf B)$ is a simple minimal free resolution of
$R/I$.

\begin{thm1} Let $I_\L=(f_1,\ldots, f_s)$ be a
lattice  ideal where $\{f_i: i=1,\ldots, s \}$ is an $R$-sequence
of binomials such that ${\bf b}_i=\deg_\A(f_i)$ are incomparable.
Let $(\bf K_\bullet,\bf \phi)$ be the Koszul complex on the $f_i$
and let $\bf B$ be the canonical system of bases  of $\bf K$. Then
$\bf (K_\bullet,\bf \phi, B)$ is an indispensable complex of
$R/I_\L$.
\end{thm1}

\begin{proof}
Let $f_i=x^{{\bf u}_i}-x^{{\bf v}_i}$. We note that if $e_J\in
B_t$ then \[\deg_{\A}(e_J)=\sum_{i\in J} \deg_A f_i\] and $(\bf
K_\bullet,\bf \phi)$ is $\A$-homogeneous. The incomparability
assumption on the degrees of the $f_i$ shows that each $\bf b_i$
is minimal and that  $\b_{1,{\bf b}_i}(R/I)=1$. It follows that
$f_i$ is an indispensable binomial, see \cite[Corollary
3.8]{ChKT}. We also note that for each $i$, $C_{\bf b_i}$ consists
of exactly two monomials.

Let $({\bf G_\bullet, \th, W})$ be a simple minimal resolution of
$R/I_\L$ where ${\bf W}=(W_0,\ldots, W_s)$ and $W_0=\{1\}$. We let
$\om_0=\id_R: K_0\mapright{} G_0$. We prove that there is a based
isomorphism $\om: \bf (K_\bullet,\bf \phi, B)\mapright{}$ $({\bf
G_\bullet, \th, W})$ which extends $\om_0$ by showing that if
$\om_i: K_i\mapright{} G_i$ has been defined for $i\le k$ then
$\om_{k+1}$ can be constructed with the desired properties. Thus
we assume that for each basis element $e_J$ of $B_k$   there
exists $c_J\in \Bbbk^*$ and $H_J\in W_k$ such that $\om_k(e_J)=c_J
H_J$.
 We note that if $e_L\in B_{k+1}$ then
 $\om_k\phi_{k+1}(e_L)$,i.e.
\[ \sum_{j\in L} \sgn[j,L]\ f_{j}\ c_{L_j}\
 H_{L_j} \]
is a simple $k$-syzygy with respect to $W_k$. This follows as in
the proof of  \cite[Corollary 3.8]{ChKT}. We will define
$\om_{k+1}: K_{k+1}\mapright{} G_{k+1}$ by specifying its image in
the basis elements $e_L$ of $B_k$ so  that the following identity
holds:
\[ \th_{k+1}\om_{k+1}(e_L)=\om_{k}\phi_{k+1}(e_L)\ .\]

\noindent Since $\om_{k}\phi_{k+1}(e_L)$ is a $k$-syzygy, it
follows that
\[ \om_{k}\phi_{k+1}(e_L)=\sum_{i=1}^{t} \th_{k+1}(p_i H_i)\]
where $H_i\in W_{k+1}$ and $\deg_{\A}(p_iH_i)=\deg_A(e_L)$. We
will show that $t=1$. First we notice that for some $i$
\[ S_{W_k}( \th_{k+1}(p_i H_i)) \cap S_{W_k}(\om_k
(\phi_{k+1}(e_L)))\neq\emptyset\ .\] Without loss of generality we
can assume that this is the case for $i=1$ and we  write $H$ in
place of $H_1$. Moreover we can assume that
\begin{itemize}
\item{} $L=\{1,\ldots, k+1 \}$ and that

\item{} $x^{{\bf u}_1} H_{L_1}\in
S_{W_k}( \th_{k+1}(H)) \cap S_{W_k}(\om_k (\phi_{k+1}(e_L)))$.

\end{itemize}

\noindent Let $q_{L_1}$ be the coefficient of
 $H_{L_1}$ in $\th_{k+1}(H)$. We have that $\deg_\A(p_1 q_{L_1})=\bf
 b_1$. We will show that $p_1 q_{L_1}$ is a constant multiple of $f_1$.
 For $t\in
L_1$ we write $L_{1,t}$ for the set $L_1\setminus \{t\}$. Since
$\th_k\th_{k+1}(H)=0$ the coefficient of $H_{L_{1,t}}$ in
$\th_k\th_{k+1}(H)$ must be zero for any $t\in L_1$. The
contributions to this coefficient  come
 from the differentiation of the term of
 $\th_{k+1}(H)$ involving $H_{L_1}$ and all other  terms of
 $\th_{k+1}(H)$ involving
 $H_{L'}$  where
 $L'\setminus \{ t'\}=L_{1,t}$.
Let $X$ be the set consisting of such   $L'$ and let $q'$ be the
coefficient of $H_{L'}$  when $L'\in X$. We get
\[0=\sgn[t, {L_1}] q_{L_1} f_t +
\sum_{L'\in X} q' \sgn[t', L'] f_{t'} \ .\] \noindent Since
$f_1,\ldots, f_s$ is a complete intersection it follows that
$q_{L_1}\in$ $\langle f_{t'}:$ $\ L'\in X\rangle$. Therefore ${\bf
b}_1\geq \deg_A(q_{L_1})\geq \deg_A(f_{t'})$ for at least one $
{t'}$. By the incomparability of the degrees of the $f_i$ it
follows that $t'=1$,  and $q_{L_1}$ is a constant multiple of
$f_1$ and thus $p_1\in \Bbbk^*$. Moreover we have shown that for
each $t$ in $L_1$ there is a term in $\th_{k+1}(H)$ involving
$H_{\{1,\ldots, \hat t, \ldots, k+1\}}$. By a degree consideration
it follows that the coefficient of this term has degree $\bf b_t$
and thus repeating the above steps we can conclude that the
coefficient of this term is a constant multiple of $f_t$. It
follows that
\[S_{W_k}(\om_k\phi_{k+1}(e_L))\subset S_{W_k}(\th_{k+1}(H))\ .\]
Since $\th_{k+1}(H))$ is simple it follows that
\[S_{W_k}(\om_k\phi_{k+1}(e_L))= S_{W_k}(\th_{k+1}(H))\ ,\] and
$\th_{k+1}(H))=c \ \om_k\phi_{k+1}(e_L)$ where $c\in \Bbbk^*$. We
let $H_L=H$ and $c_L=c^{-1}$. It follows that the homomorphism
$\om_{k+1}:K_{k+1}\mapright{} G_{k+1}$ defined  by setting
\[\om_{k+1}(e_L)=c_L H_L\ \]
has the desired properties.
\end{proof}

Next we consider strongly indispensable complexes.

\begin{def1} \label{def_stongly_ind_complex}{\rm Let $\bf(F_\bullet, \phi,B)$ be a based complex.
We say that $(\bf F_\bullet, \phi, B)$  is a {\it  strongly indispensable
complex} of $R/I$ if
 for every based minimal   free resolution $\bf (G_\bullet,\theta,C)$ of
       $R/I$, (not necessary simple) with $C_0=\{1\}$, there is an injective based homomorphism
        $\bf \omega: \bf(F_\bullet, \phi,B)
       \mapright{}
       \bf (G_\bullet,\theta,C)$
       such that $\omega_0 = id_R$.
If ${\bf B}=(B_j)$ and $E\in B_{i+1}$ we say that
$\phi_{i+1}(E)\in F_i$ is a \emph{strongly indispensable i-syzygy}
of $R/I$.}
 \end{def1}

Strongly indispensable complexes are indispensable. This is a
strict inclusion as \cite[Example 6.5]{ChTh} shows. When $I_\L$ is
a lattice ideal,  the algebraic Scarf complex \cite[Construction
3.1]{PS1}, is shown to be ``contained" in the minimal free
resolution of $R/I_\L$, \cite[Theorem 3.2]{PS1}, and is a strongly
indispensable complex. Moreover as follows from Corollary
\ref{start_indispens_complex} the strongly indispensable
$0$-syzygies of $R/I_\L$ coincide with the indispensable
$0$-syzygies of $R/I_\L$ and are   the indispensable binomials of
$I_\L$. For higher homological degrees this is no longer the case.
First we note the following:

\begin{thm1}\label{strong_indisp} Let $I_\L$ be a lattice ideal
and let $(\bf F_\bullet, \phi, B)$ be a strongly indispensable
complex for $R/I_\L$.  Let ${\bf B}=(B_j)$, $E\in B_{i+1}$ and
$\deg_\A(E)=\bf b$. Then $\dim_{\Bbbk} \HH_i(\D_{\gcd}({\bf
b}))=1$ and $\bf b$ is a minimal $i$-Betti degree of $R/I_\L$.
\end{thm1}

\begin{proof} Suppose that  $\dim_{\Bbbk}
\HH_i(\D_{\gcd}({\bf b}))>1$ or that there is an $i$-Betti degree
$\bf b'$ such that $\bf b'< \bf b$. Let $({\bf G_\bullet, \th, C})$ be a
minimal resolution of $R/I_\L$ where ${\bf C}=(C_i)$, let $\bf
\omega: \bf(F_\bullet, \phi,B) \mapright{}
       \bf (G_\bullet,\theta,C)$ be the based homomorphism of Definition
\ref{def_stongly_ind_complex} and suppose that $\om(E)=c H$ where
$H\in C_{i+1}$ and $c\in \Bbbk^*$. By our assumptions there exists
$H'\in C_{i+1}$ such that $H'\neq H$ and $\deg_\A(H)\le \bf b$.
Let $x^{\bf a}\in C_{ \bf b-\bf b'}$. By replacing $H$ with
$H+x^{\bf a} H'$ we get a new basis $C'_{i+1}$ of $G_{i+1}$ and a
new system of bases $\bf C'=(C_j')$, where $C'_j=C_j$ for $j\neq
i+1$. Let $\om': (\bf F_\bullet,\phi, B)\mapright{} \bf (\bf
G_\bullet,\th, C')$ be the based homomorphism of Definition
\ref{def_stongly_ind_complex}: $\om_j=\om_j'$ for $j\le i$. Let
$H''\in C_{i+1}'$ be such that  $\om_{i+1}'(E)=c'H''$ where $c'\in
\Bbbk^*$. Thus
\[ \th_{i+1}(c' H'')=\th_{i+1}(\om_{i+1}'(E))=
\om_{i}'\phi_{i+1}(E)=\om_{i}\phi_{i+1}(E)=\]
\[\th_{i+1}(\om_{i+1}(E))=c\th_{i+1}(cH)\ . \] It follows that
$c'H''-cH \in \ker \th_{i+1}$. If $H''\neq H+x^{\bf a} H'$ then
$H''\in C_{i+1}$ and we get a direct contradiction to  the
minimality of $(\bf G_\bullet,\th,C)$. If $H''= H+x^{\bf a} H'$
then $\th_{i+1}((c'-c)H+c'x^{\bf a} H')=0$. Examination of the two
cases when (a)
 $c'\neq c$, and (b)
 $c'=c$,
leads again to a contradiction of the minimality of the resolution
$(\bf G_\bullet,\th,C)$.
\end{proof}

Theorem \ref{strong_indisp} shows that the two conditions
\begin{enumerate}
\item{}$\dim_{\Bbbk} \HH_i(\D_{\gcd}({\bf b}))=1$
\item{}$\bf b$ is a minimal $i$-Betti degree
\end{enumerate}
are necessary for the existence of a strongly indispensable
$i$-syzygy in $\A$-degree $\bf b$. The following example shows
that these conditions are not sufficient for the existence of an
indispensable $i$-syzygy and consequently of a strongly
indispensable $i$-syzygy.

\begin{example}{\rm Consider the lattice ideal $I_\L= \langle f_1,
f_2 \rangle$ where $f_1= x_1-x_2$, $f_2=x_2-x_3$ and $\deg_\A
f_i=1$. Let $(\bf K_\bullet, \phi)$ be the Koszul complex on the $f_i$. By
considering the $i$-Betti numbers for $i=1, 2$ it is immediate
that $\dim_{\Bbbk} H_2(\D_2)=1$ and $2$ is a minimal $2$-Betti
degree. However there is no indispensable complex of length
greater than 0, since the generators of $I_\L$ are not indispensable
binomials.}
\end{example}

Generic lattice ideals are characterized by the condition that the
binomials in a minimal generating set have full support,
\cite{PS1}. In this case the Scarf complex is a minimal free
resolution of $R/I_\L$  and each of the Betti degrees of $R/I_\L$
satisfy the conditions of Theorem \ref{strong_indisp}. We finish
this section by giving the strongest result for the opposite
direction of Theorem \ref{strong_indisp}.

\begin{thm1} \label{strongly} Let $I_\L$ be a lattice ideal. The  $\A$-homogeneous
minimal free resolution $(\bf F_\bullet, \phi, B)$ of $R/I_\L$ is
strongly indispensable
 if and only if for each $i$-Betti degree $\bf b$ of
 $R/I_\L$,
$\bf b$ is a minimal $i$-Betti degree  and $\dim_{\Bbbk}
\HH_i(\D_{\gcd}({\bf b}))=1$.
\end{thm1}

\begin{proof} One direction of this theorem follows directly from
Theorem  \ref{strong_indisp}. For the other direction we assume
that $\bf b$ is minimal whenever $\bf b$ is an $i$-Betti degree
and that $\dim_{\Bbbk} \HH_i(\D_{\gcd}({\bf b}))=1$ for all $i$.
 Let $({\bf G_\bullet, \th, D})$ be a
minimal free resolution of $R/I_\L$. By assumption the
$\A$-degrees of the elements of $D_i$ are distinct and
incomparable. It follows that the $\A$-homogeneous isomorphism
$\om: \bf F\mapright{} \bf G$ that extends $\id_R: F_0\mapright{}
G_0$ is a based homomorphism.
\end{proof}

{\bf Acknowledgment}

The authors would like to thank Ezra Miller for his essential
comments on this manuscript.


\begin{thebibliography}{00}

\bibitem{ATY} S. Aoki,  A. Takemura and R. Yoshida,  {\em
Indispensable monomials of toric ideals and Markov bases}, Journal
of Symbolic Computation \textbf{43} (2008) 490-507.

\bibitem{AH} A. Aramova and J. Herzog, {\em Koszul cycles and
Eliahou-Kervaire type resolution}, J. Algebra  \textbf{181} (1996)
347-370.




\bibitem{BaSt} D. Bayer, B. Sturmfels, {\em Cellular resolutions
of monomial modules}, J. Reine Angew. Math. \textbf{502} (1998),
123-140.



\bibitem{BCMP2} E. Briales, A. Campillo, C. Marijuan and P. Pis$\acute{o}$n,
{\em Combinatorics of syzygies for semigroup algebra}, Collectanea
Mathematica \textbf{49} (1998) 239-256.

\bibitem{BH} W. Bruns and J. Herzog, {\em Semigroup rings and
simplicial complexes}, J. Pure Appl. Algebra \textbf{122} (1997)
185-208.

\bibitem{BW} A.~Bjorner and J.~W.~Walker {\em A homotopy complementation formula
for partially ordered sets}, European J. Combin. \textbf{4}
(1983), 11--19.

\bibitem{CG} A. Campillo and Ph. Gimenez, {\em Syzygies of affine toric varieties}, J. Algebra 225 (2000) 142-161.

\bibitem{CP} A. Campillo and P. Pis\`{o}n,  {\em L'id{\'e}al d'un
semi-group de type fini},  Comptes Rendues Acad. Sci. Paris,
S{\'e}rie I,  \textbf{316} (1993) 1303-1306.


\bibitem{ChKT} H. Charalambous, A. Katsabekis, A. Thoma,  {\em Minimal systems of binomial
generators and the indispensable complex of a toric ideal},  Proc.
Amer. Math. Soc. \textbf{135} (2007) 3443-3451.

\bibitem{ChTh} H. Charalambous, A. Thoma,  {\em On the genearalized Scarf complex for lattice ideals},
preprint.



\bibitem{DS} P. Diaconis and B. Sturmfels,  {\em Algebraic
algorithms for sampling from conditional distributions},  Ann.
Statist.,  {\bf 26} (1) (1998) 363-397.

\bibitem{ES} D. Eisenbud and B. Sturmfels, {\em Binomial ideals}, Duke Math. J. {\bf 84} (1996) 1-45.



\bibitem{MS} E. Miller and B. Sturmfels,  Combinatorial Commutative
Algebra, Graduate Texts in Mathematics \textbf{227}  Springer
Verlag, New York 2005.

\bibitem{Hi-O} H. Ohsugi and T. Hibi,  {\em Indispensable binomials
of finite graphs},  J. Algebra Appl. {\bf 4} (2005),  no 4,
421-434.

\bibitem{H-O} H. Ohsugi and T. Hibi,  {\em Toric ideals arising from contingency
tables}, Proceedings of the Ramanujan Mathematical Society's
Lecture Notes Series,  (2006) 87-111.

\bibitem{OjVi} I. Ojeda and A. Vigneron-Tenorio,
{\em Simplicial complexes and minimal free resolution of monomial
algebras}, preprint, arXiv:0810.4836v1

\bibitem{P} P. Pis\`{o}n Casares,  {\em The short resolution of a lattice ideal}, Proc. Amer. Math. Soc. 131 (2003)
1081-1091.

\bibitem{PVT} P. Pis\`{o}n Casares and A. Vigneron-Tenorio,
{\em First syzygies of toric varieties and diophantine equations
in congruence}, Commun. Algebra \textbf{29} (2001) 1445-1466.


\bibitem{PS1} I. Peeva and B. Sturmfels,  {\em Generic lattice
ideals},  J. Amer. Math. Soc. \textbf{11} (1998) 363-373.

\bibitem{PS2}  I. Peeva and B. Sturmfels,  {\em Syzygies of
codimension 2 lattice ideals},  Math Z. \textbf{229} (1998) no 1,
163-194.

\bibitem{Ro} J. Rotman,  An introduction to Algebraic
Topology, Graduate Texts in Mathematics \textbf{119}  Springer
Verlag, New York 1988.

\bibitem{Stanley} R. Stanley, Combinatorics and commutative
algebra, Progress in Mathematics \textbf{41}, Birkh\"{a}user,
Boston 1996.

\bibitem{St} B. Sturmfels,  Gr{\"o}bner Bases and Convex Polytopes.
University Lecture Series,  No. 8 American Mathematical Society
Providence,  R.I. 1995.

\bibitem{AT}  A. Takemura and S. Aoki,  {\em Some characterizations
of minimal Markov basis for sampling from discrete conditional
distributions},  Ann. Inst. Statist. Math.,  {\bf 56} (1)(2004)
1-17.
\end{thebibliography}
\end{document}